\documentclass[12pt]{article}
\usepackage{amsthm,amsfonts,amsmath,amssymb}
\usepackage{graphicx}
\usepackage{ tabls}
\usepackage{marvosym}

\usepackage{pstricks,pst-node,pst-tree, pstcol}

\newcommand{\CC}{\boldsymbol C}
\newcommand{\PP}{\boldsymbol P}
\newcommand{\QQ}{\boldsymbol Q}

\newcommand{\Cc}{\boldsymbol c}

\newcommand{\Pp}{\boldsymbol p}
\newcommand{\Qq}{\boldsymbol q}

\newcommand{\HH}{\boldsymbol M}

\newcommand{\s}{{\large s}}
\newcommand{\m}{{\large h}}

\newcommand{\ts}{{\tt s}}
\newcommand{\tm}{{\tt h}}

\newcommand{\swsw}{{\tt ss}}
\newcommand{\masw}{{\tt hs}}
\newcommand{\swma}{{\tt sh}}
\newcommand{\mama}{{\tt hh}}

\begin{document}
\title{
The Monty Hall Problem\\
 in the Game Theory Class} 
\author {Sasha Gnedin\thanks{A.V.Gnedin@uu.nl} } 
\maketitle

\section{Introduction}

\begin{itemize}
\item[]
{\it Suppose you're on a game show, and you're given the choice of three doors: 
Behind one door is a car; behind the others, goats. 
You pick a door, say No. 1, and the host, who knows what's behind the doors, 
opens another door, say No. 3, which has a goat. He then says to you, 
``Do you want to pick door No. 2?" Is it to your advantage to switch your choice?}
\end{itemize}
 With these famous words the {\it Parade Magazine} columnist vos Savant opened an exciting chapter in mathematical didactics. 
The puzzle, which has become known as the {\it Monty Hall Problem} (MHP), has broken the records of popularity
among the probability paradoxes.
The book by Rosenhouse \cite{Rosenhouse} and the   Wikipedia entry on the MHP present the history  and variations of the problem.

In the basic version of the MHP the rules of the game specify that the host must always reveal one of the
unchosen doors to show that there is no prize there. Two remaining unrevealed doors  
may hide the prize, creating  the illusion of symmetry and suggesting that the action does not matter.
 However, the symmetry is fallacious, and 
 switching is a better action, doubling the probability of winning.

There are two main explanations of the paradox. One of them, simplistic, amounts to just counting 
the mutually exclusive cases: either you win with switching or with holding the first choice.
A more sophisticated argument, included in the textbooks as an exercise on the Bayes theorem, calculates the conditional
probability of winning in the situation described in vos Savant's wording of the problem.

In \cite{Gill} a critical analysis of the conventional approaches to the MHP has been done, with the advocated viewpoint 
that the whole situation of making decision, viewed as a multistep process, 
is a challenging instance of mathematical modelling, very much amenable to the analysis within the game-theoretic framework.
The textbook by Haggstr\"om \cite{Haggstrom}  puts the zero-sum game in matrix form and presents a minimax solution. 
Further steps in this direction were done in \cite{Doors} and \cite{Dominance}, where it was argued that
the game-theoretic concept of dominance allows to analyse the problem under fairly general assumptions  on the prior 
information of the decision-maker, including the very interesting case of nouniform distribution, only  
occasionally included in exercise sections
of probability textbooks \cite{Grinstead, Tijms}.

One elementary new observation we make is that for every contestant's strategy of playing the game
after the initial door has been chosen, there is always (at least)
one `unlucky'  door, the same for every  admissible algorithm for revealing the doors by the host.
The prize is never found behind the `unlucky' door, hence
the contestant will always lose in at least one case out of three.
This leads, rather straightforwardly,
to the worst-case winning probability 2/3. 
We prefer, however, for the sake of instruction and in anticipation of future 
generalizations to give here a full Bayesian analysis, with the host biased in any possible way.

These notes, expanding upon the cited literature are intended to show that the MHP is indeed
an excellent occasion to expose the undergraduate students to  basic ideas of the game theory and
to decision-making under various uncertainty scenarios. 
Due to its remarkable  symmetry features  the zero-sum version of the game will  undoubtedly enter the Hall of Fame of the classical  games 
like {\it matching pennies} and {\it paper-scissors-stone}.

\section{The MHP as a sequential decision process}

The game involves two actors which we christen Monte and Conie.
In the basic  scenario to follow, the prize is hidden behind  one of three doors by Monte, 
then 
Conie picks a door which is kept unrevealed,
then one of the unchosen doors is revealed as not containing the prize, and then an offer is made to switch the choice
from the initial guess to another unrevealed door. Conie  wants the prize 
and she wins it if her final choice falls on the door where the prize was hidden.

To state the rules  more rigorously and to introduce some notation for the admissible actions we number the doors 1,2,3. 
The game  consists of four moves:
\begin{itemize}
\item[(i)]  Monte 
chooses a door number $\theta$ from doors 1,2,3 to hide the prize. He keeps $\theta$ in secret.
\item[(ii)] Conie picks door number $x$ out of 1,2,3. Monte observes Conie's choice.
\item[(iii)] Monte reveals a door distinct from $x$ as not hiding the prize,
and offers to Conie the possibility of choosing between door $x$ and another unrevealed door $y$.
\item[(iv)] Conie finally chooses door $z$ from doors $x$ and $y$.  Conie wins the prize if $z=\theta$, otherwise she wins nothing.
\end{itemize}

Conie's ignorance about the location of the prize  means in the mathematical language that her actions cannot depend  on $\theta$ explicitly.  
There is also another, more subtle, indefinite factor of which Conie is ignorant:
this is the way Monte chooses between two doors to reveal
in the case $x=\theta$. These two indefinite factors, which are under the control of her opponent, comprise
Conie's decision-making environment.

We say that there is a match if $\theta=x$, in which case $y$ on step
(iii) 
 can be one of two numbers distinct from $\theta$; whereas if there is a mismatch,
$\theta\neq x$, the rules force $y=\theta$. On step (iv), we say that Conie takes action `hold' 
(of sticking with the initial choice), denoted $\tm$, if $z=x$; 
and that she takes action `switch' (from the initial choice), denoted $\ts$, if $z=y$.

\par To illustrate, a sequence of admissible moves $\theta,x,y,z$ could be 2212, which means that
Monte hides the prize behind door 2, Conie picks door 2 (so there is a match),
then Monte offers to switch to door 1 (by revealing door 3 as not containing the prize), and finally 
Conie plays $\tm$ by sticking with her initial choice 2. Since $z=\theta~$ Conie wins the prize.

\begin{figure}
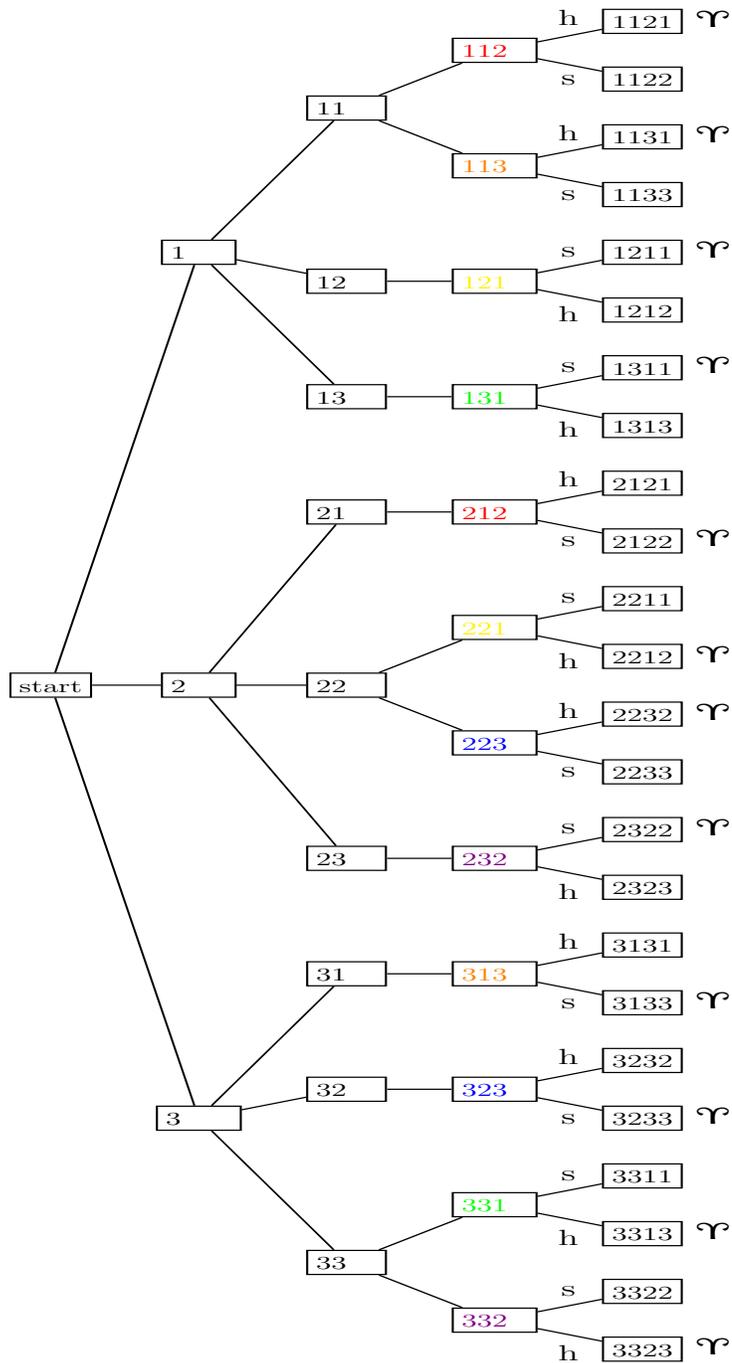

$
{
\resizebox*{10cm}{18cm}
{
\newcommand{\XX}[1]{%
\Tr{\psframebox{\rule{0pt}{9pt}#1}}}
\psset{angleB=90,angleA=90}
\pstree
[treemode=R]
{\XX{start}}
{
\pstree{{\XX{1~~~~}}}
{
\pstree{{\XX{11~~~}}}{{\pstree{{\XX{\color{red}112~~}}}{{\XX{1121}~{\large\Aries}}\taput{\m}{\XX{1122}}}}\tbput{\s}
\pstree{{\XX{\color{orange}113~~}}} {{\XX{1131}~{\large\Aries}}\taput{\m}{\XX{1133}}}}\tbput{\s}
\pstree{{\XX{12~~~}}}{\pstree{{\XX{\color{yellow}121~~}}}{{\XX{1211}~{\large\Aries}}\taput{\s} {\XX{1212}}}}\tbput{\m}
\pstree{{\XX{13~~~}}}{\pstree{{\XX{\color{green}131~~}}}{{\XX{1311}~{\large\Aries}}\taput{\s}{\XX{1313}}}}\tbput{\m}
}

\pstree{{\XX{2~~~~}}}
{
\pstree{{\XX{21~~~}}}{\pstree{{\XX{\color{red}212~~}}}{{\XX{2121}}\taput{\m}{\XX{2122}~{\large\Aries}}}}\tbput{\s} 
\pstree{{\XX{22~~~}}}
{\pstree{{\XX{\color{yellow}221~~}}}{{\XX{2211}}\taput{\s} {\XX{2212}~{\large\Aries}}}\tbput{\m} \pstree{{\XX{\color{blue}223~~}}}
{{\XX{2232}~{\large\Aries}}\taput{\m} {\XX{2233}}}}\tbput{\s} 
\pstree{{\XX{23~~~}}}{\pstree{{\XX{\color{violet}232~~}}}{{\XX{2322}~{\large\Aries}}\taput{\s} {\XX{2323}}}}\tbput{\m} 
}
\pstree{{\XX{3~~~~ }}}
{
\pstree{{\XX{31~~~}}}{\pstree{{\XX{\color{orange}313~~}}}{{\XX{3131}}\taput{\m} {\XX{3133}~{\large\Aries}}}} \tbput{\s}
\pstree{{\XX{32~~~}}}{\pstree{{\XX{\color{blue}323~~}}}{{\XX{3232}}\taput{\m}  {\XX{3233}~{\large\Aries}}}} \tbput{\s}
\pstree{{\XX{33~~~}}}{\pstree{{\XX{\color{green}331~~}}}{{\XX{3311}}\taput{\s}   {\XX{3313}~{\large\Aries}}}\tbput{\m}
\pstree{{\XX{\color{violet}332~~}}}{{\XX{3322}}\taput{\s}{\XX{3323}~{\large\Aries}}}}\tbput{\m} 
}
}
}
}
$
\caption{The game tree with succession of moves $\theta,x,y,z$, and Conie's winning terminal positions marked \Aries.}
\label{TheBigGameTree}
\end{figure}

Positions in the game represent all substantial information available for the move to follow.
These are represented by vertices in the game tree in Figure \ref{TheBigGameTree}. 
An edge connects a position to another position achievable in one move.
The play starts at the root vertex with Monte's move leading to a position $\theta\in\{1,2,3\}$, then the moves of actors alternate. 
A path in the tree from the root to a terminal vertex is determined by the actions of {\it both} actors.
Each path directed away from the root ends in a leaf node, with Conie's winning positions $(\theta,x,y,z)$ being those 
with $z=\theta$.

There is one important feature of the game which we indicate by coloring positions in the figure. 
Conie does not know the 
winning door $\theta$. The information of Conie on her second  move can be specified by partitioning the collection of relevant positions
in {\it information sets}.  For the first Conie's move $x$ there is just one information set $\{1,2,3\}$.
On her second move,
 Conie cannot distinguish e.g. between  positions 121 and 221, since for $x=2$ and $y=1$
(when the revealed door is 3) the prize can be  behind any of the doors 1 and 2; thus $\{121, 221\}$ is one information set
which we denote $*21$, with $*$ staying for the unknown admissible value of $\theta$ (1 or 2 in this case).

The complete list of information sets for
the second Conie's move is $*12, *13, *21, *23, *31, *32$.
Possible  moves are labelled by actions $\{\tm,\ts\}$, and the action depends on the set. Thus if the action is $\tm$ from 
position $112$, then it must be $\tm$ from position $212$.
The game tree with partition of positions into information sets is sometimes called the {\it Kuhn tree}.

Monte always knows the current position in the game exactly,
so his information sets are singletons.

\section{Strategies and the payoff matrix}

A  {\it strategy} of Monte is a rule which for each position, where Monte is on the move, specifies a follower exactly.
For his first move from the starting position a strategy specifies a value of $\theta$. 
For his second move a strategy specifies the value of $y$  for each  position $(\theta,x)$,
so we can consider $y$ as a function $y= d_\theta(x)$, where
\begin{eqnarray*}
 d_\theta(x)=\theta ~~~{\rm if ~~~} x\neq\theta,\\
 d_\theta(x)\neq \theta ~~~{\rm if~~~} \theta=x.
\end{eqnarray*}
Simply put, Monte's strategy can be encoded in a pair of door numbers like $12$, where $\theta=1$ is the door hiding the prize, and 
$d_1(1)=2$ is the door to which the switch is offered in the case of match. This indeed determines the function $d_1(\cdot)$ uniquely because
$d_1(2)=d_1(3)=1$. With this notation,
the complete list of Monte's strategies is 

$${\mathcal M}= \{12,~ 13,~ 21,~ 22,~ 31,~ 32\}.$$

A {\it strategy} of Conie is a rule which for each position where Conie is on the move 
specifies a follower, in a way consistent with partition in information sets. 
Thus,  her strategy must specify the value of $x$. Furthermore, for each initial choice $x$ and the door offered for switch $y\neq x$ 
her strategy must specify an action from the set $\{\tm,\ts\}$, which  
is a function $a_x(\cdot)$, so that the second Conie's move is    
\begin{eqnarray*}
 z=x ~~~{\rm if ~~~} a_x(y)=\tm,\\
  z=y ~~~{\rm if ~~~} a_x(y)=\ts.
\end{eqnarray*}
When $x$ is fixed $a_x(y)$ must be specified for two possible values of $y$. We can write therefore Conie's strategy as a triple
like $2\ts\tm$ which specifies $x=2$ and $a_2(1)=\ts, a_2(3)=\tm$.
The second entry $\ts$ of $2\ts\tm$ encodes the action for smaller door number $y$, and the third entry $\tm$ for larger.
With similar conventions, the complete list of twelve strategies of Conie is \\
$${\mathcal C}= \{1\ts\ts,~ 1\ts\tm, ~1\tm\ts, ~1\tm\tm, ~2\ts\ts, ~2\ts\tm, ~2\tm\ts, ~2\tm\tm, ~3\ts\ts, ~3\ts\tm, ~3\tm\ts, ~3\tm\tm\}.$$

\noindent
Every strategy of the kind $x\tm\tm$ or $x\ts\ts$ will be called {\it constant-action} strategy,
and every strategy of the kind $x\ts\ts$ will be called {\it always-switching} stratgey.
Strategies $x\ts\tm$ and $x\tm\ts$ are {\it context-dependent} strategies.

A strategy profile of the actors is a strategy of Monte and a strategy of Conie. When a strategy profile is fixed,
the course of the game, represented by a path in the game tree, is fully inambiguous.
For instance, when the profile (12,2\ts\tm) is played by the actors, $\theta,x,y,z$ is
1211, because Monte offers a switch to door 1,  and since 1 is the smaller door number of possible values of $y\in\{1,3\}$ Conie reacts with $\ts$,
hence 
winning the prize. The second entry 2 of Monte's strategy 12 is immaterial for this outcome since there is a mismatch $\theta\neq x$.
We adopt the convention that Conie receives payoff  1 when she wins the prize and 0 otherwise. 
The matrix $\CC$ in Figure \ref{ConPayoff} 
represents the correspondence between the strategy profiles and Conie's payoffs.

\begin{figure}
\begin{center}
\begin{tabular}{c|cccccc}
$\theta,y=$    & 12 & 13 & 21 & 23 & 31 & 32\\

\hline
1\swsw &0  & 0 &  1  &  1 & 1  &1\\
1\masw &1  & 0 &  0 &  0 & 1  &1\\
1\swma &0  & 1 &  1  &  1 & 0  &0\\
1\mama &1  & 1 &  0  &  0 & 0  &0\\
    &   &   &     &    &    &  \\
2\swsw &1  & 1 &  0  &  0 & 1  &1\\
2\masw &0  & 0 &  1 &  0 & 1  &1\\
2\swma &1  & 1 &  0  &  1 & 0  &0\\
2\mama &0  & 0 &  1  &  1 & 0  &0\\
    &   &   &     &    &    &  \\
3\swsw &1  & 1 &  1  &  1 & 0  &0\\
3\masw &0  & 0 &  1 &  1 & 1  &0\\
3\swma &1  & 1 &  0  &  0 & 0  &1\\
3\mama &0  & 0 &  0  &  0 & 1  &1\\
\end{tabular}
\end{center}
\caption{Conie's payoff matrix $\CC$} 
\label{ConPayoff}
\end{figure}

Unlike the game tree, the simplified matrix representation ignores many substantial
attributes like moves,  positions and information sets. The game in this matrix form amounts to a very simple procedure:
Monte picks a column and Conie picks a row of matrix $\CC$. Conie's payoff is then the entry of $\CC$ that stays
on the intersection of the selected row and column.

An advantage of the matrix representation is that it unifies and simplifies
the analysis. In particular, we can compare different Conie's strategies under different assumptions on the decision-making 
environments, that is Monte's behaviors. For instance, we can think of two Conies which simultaneously 
play, say, 1\masw~ and 3\swsw~ (respectively) against the same strategy of Monte.
 Note that `simultaneous' in this context refers to  a logical comparison of the
outcomes, rather than to a real play.

A quick inspection of the matrix $\CC$ shows that the problem has a property called {\it weak dominance}: 
\begin{itemize}
\item[]
{ for every Conie's strategy $A$ which in some situation employs action $\tm$ (so $A$  is of the kind $x\tm\tm$, $x\tm\ts$ or $x\ts\tm$) there exists
an always-switching strategy $B$ such that if $A$ wins against  $S$, 
then $B$ wins against $S$ as well, whichever strategy of Monte $S$}.
\end{itemize}
For instance, strategy 1\tm\ts~ which does not switch to door $y=2$, is dominated by always-switching strategy 2\ts\ts.
Amazingly, $1\tm\ts$ is beaten in the situation $\theta=1, y=3$ where the strategy uses the {\it switch} action!
The never-switching strategy 1\tm\tm~ is dominated by both 2\ts\ts~and~  3\ts\ts.
The dominance is a very strong argument for excluding the strategies which may employ the action $\tm$.

The dominance feature has a nice interpretation in terms of `unlucky' door.
Suppose for a time-being that the prize-hiding is a move of nature (or quiz-team), out of control of Monte.
All Monte can do is to choose door $y$ to reveal when he can.
\begin{itemize}
\item[]
{The `unlucky' door theorem says: whichever Conie's strategy $A$, there exists a door $u$ (depending on $A$) such that
under $A$ the final choice $z$ does not fall on $u$ when $\theta=u$, whichever Monte's way of revealing doors when he can.
Strategy $u\ts\ts$~then weakly dominates $A$.}
\end{itemize}

Door $u=u(A)$ is marked in $\CC$ by two zeroes in positions $u*$ of the row $A$, for instance $u=2$ for $A=1\tm\ts$, and $u=2$ or $3$ for
$A=1\tm\tm$.

So the existence of the `unlucky' door means that Conie never wins for $\theta=u$, no matter how Monty reveals the doors.
Therefore, 
\begin{itemize}
\item[]
{If $\theta$ is chosen uniformly at random Conie's winning probability cannot exceed $2/3$.}
\end{itemize}


\section{The Bayesian games}\label{Bayes}

We described above the sets of {\it pure} strategies $\mathcal M$ and $\mathcal C$, selecting from which a particular 
profile determines inambiguously a succession of moves.
A {\it mixed} strategy of an actor is an assignment of probability to each pure strategy 
of the actor.
 Playing mixed strategies means using chance devices to choose from $\mathcal M$ and $\mathcal C$. 
A mixed strategy of Monte can be written as a row  vector $\QQ$ with six components corresponding to her pure strategies.
A mixed strategy of Conie can be written as a row vector $\PP$ of length twelve.

Conie's generic mixed strategy has 11 parameters, as probabilities add to one. A general theorem due to Kuhn says that if the actors
do not forget their private history 
(perfect recall), like in our game,   it is enough to consider a smaller class of {\it behavioral} strategies.
A behavioral strategy specifies probability distribution on the set of admissible moves for every information set of an actor.
Thus behavioral strategy of Conie is described by 8 parameters, 2 numbers for distribution of $x$ and 
6 biases for coins tossed for each information set.
This subtle distinction of the general and behavioral strategies is only mentioned here for the sake of instruction,
but the distinction becomes important in some games where an actor may forget some elements of her private history.

The {\it support} of a mixed strategies are the pure strategies with nonzero probability. 
For instance, the support of  $(\dots,0,1,0,\dots)$ is one pure strategy,
thus we simply identify pure strategies with 
mixed strategies of this kind. This identification allows us to view a mixed strategy 
as a convex combination (also called mixture) of the pure strategies constituting the support.
A mixed strategy $\PP$ (respectively $\QQ$) is called {\it fully supported} if 
the support is $\cal C$ (respectively $\cal M$).

When strategy profile $(\PP,\QQ)$ is played  by the actors, the expected payoff for Conie, equal to her probability of winning the prize,  
is computed by the matrix multiplication as $\PP\CC\QQ^T$, where $^T$ denotes transposition. 
This way of computation presumes   that the actors' choices
 of pure strategies are independent random variables with values in $\mathcal M$ and $\mathcal C$ (respectively), 
which may be simulated  by the actors'	private randomization devices.
The independence of individual strategies is a feature required by the idea of {\it noncooperative} game, in which there is no communication of actors
to play a joint strategy.


In the {\it Bayesian setting} of the decision problem Monte is supposed to play some fixed strategy $\QQ$, known to Conie.
The probability textbooks consider the Bayesian setting for the  MHP, with $\QQ$ being the uniform distribution over $\mathcal M$,
which is equivalent to the assumption that Monte picks $\theta$ from $\{1,2,3\}$ by rolling a three-sided die, 
and in the event of match picks a door $y$ from two possibilities by tossing a fair coin.
The general Bayesian formulation, which may be called a game against the nature,
models the situation where Conie deals with unconcious random algorithm
which neither has own goals nor can take advantage of Conie's mistakes.
Conie's optimal 
behavior is then a {\it Bayesian strategy} $\PP'=\PP'(\QQ)$ which maximizes the probability of winning the prize:
$$\PP'\CC\QQ^T= \max_{\PP}  \PP\CC\QQ^T.$$   
Bayesian strategy always exists by 
continuity of 
$\PP\CC\QQ^T$ as function of $\PP$, and compactness of the set of mixed strategies.
Moreover, by linearity of  $\PP\CC\QQ^T$ and convexity 
\begin{itemize}
\item[]
The general 
Bayesian strategy is an arbitrary mixture of the pure Bayesian strategies.  
\end{itemize}

Now suppose $\QQ$ is fully supported, that is every pure strategy from $\cal M$ is played with positive probability,
and let $A$ and $B$ be two distinct pure strategies of Conie such that
$B$ weakly dominates $A$. Since all rows of $\cal C$ are distinct the latter means that the collection 
of columns in which row $A$ has a 1 is a proper subset of the collection of columns in which $B$ has a 1. 
But then the winning probability  is strictly larger for  $B$ than for $A$. Since every strategy admitting action $\tm$ is weakly 
dominated by some always-switching strategy we conclude that
\begin{itemize}
\item[] If $\QQ$ is fully supported then only always-switching pure strategies can be Bayesian, hence every 
Bayesian strategy is a mixture of always-switching strategies.
\end{itemize}
Thus in the `generic' case the Bayesian principle of optimality excludes strategies which may use $\tm$.
Having retained only always-switching strategies, it
is easy to find all Bayesian strategies explicitly.

\begin{itemize}
\item[]
Suppose that according to fully supported $\QQ$ the probability to hide the prize behind door $\theta\in\{1,2,3\}$ is $\pi_\theta$.
Let $\pi_1\geq\pi_2\geq\pi_3>0$ (otherwise re-label the doors).
Then the only pure Bayesian strategies are 
\begin{itemize}
\item[\rm(1)]  3\ts\ts~ if   $\pi_1\geq \pi_2>\pi_3$,
\item[\rm(2)]  2\ts\ts~ and 3\ts\ts~ if $\pi_1>\pi_2=\pi_3$,
\item[\rm(3)]    1\ts\ts,~ 2\ts\ts~ and 3\ts\ts~  if   $\pi_1=\pi_2=\pi_3=1/3$. 
\end{itemize}
\end{itemize}

If $\QQ$ is not fully supported  then 
the listed always-switching strategies are still Bayesian, but some other strategies may be Bayesian too.
For instance, a never-switching strategy $x\tm\tm$ is Bayesian if (and only if) Monte always hides the prize behind door $x$, although
even then every always-switching strategy $x'\ts\ts$ with $x'\neq x$ is Bayesian also.

More interestingly, suppose $\QQ=(1/3,0,1/3,0,1/3,0)$. This is a model for a `crawl' behavior  of the host \cite{Rosenthal}, who is eager
to reveal the door with higher number when he has a choice. In this case the Bayesian pure strategies are
$1\ts\ts,~ 2\ts\ts,~3\ts\ts$ and $1\tm\ts,~ 2\tm\ts,~3\tm\ts$.
In general, the rule to determine all Bayesian pure strategies by excluding the dominated strategies is the following:
\begin{itemize}
\item[]
If strategy $\theta,y$ with $y=\min\{1,2,3\}\setminus\{\theta\}$ enters
$\QQ$ with nonzero probability then $\theta\ts\tm$ is not Bayesian.

\item[] If strategy $\theta,y$ with $y=\max\{1,2,3\}\setminus\{x\}$ enters $\QQ$ with nonzero probability then  $\theta\tm\ts$ is not Bayesian.
\end{itemize}

For arbitrary $\QQ$~ Conie's strategy `point at the door which is the least likely to hide the prize, then always switch' is Bayesian, 
giving the probability of win $1-\min_{\theta} \pi_{\theta}$, where $\pi_{\theta}$ is the probability of $\theta=1,2,3$. 
Remarkably, all what Conie needs to know to play optimally in a Bayesian game is the number of the least likely door. 
She tips at this door and switches all the time, winning the complementary probability no matter which are the biases  
for revealing the doors in case of match.

In the special case $\pi_1=\pi_2=\pi_3=1/3$ every always-switching strategy is Bayesian no matter how 
$y$ is chosen in the event of match $x=\theta$; a conclusion usually shown in the literature with the 
help of conditional probabilities.

In fact, for an optimal choice of $x$, the {\it conditional} probability of winning with switching in every position
$xy$ is at least $1/2$. This is implied by the overall optimality (in the Bayesian sense), and is a simple instance of the
general Bellman's dynamic programming principle.

\section{The zero-sum game}

The zero-sum game is a model for interaction of actors with antagonistic goals. 
Conie wins the prize when Monte loses it. 
The essense of the zero-sum game is the {\it worst-case analysis}.
What Conie can guarantee in the game, when the behavior of Monte can be arbitrary (but agreeing with the game rules)?
What is the worst behavior of Monte?

We can write Monte's payoffs as another $12\times 6$ matrix $\HH$, but this is not necessary as $\HH=-\CC$, so 
$\CC$ contains all information about the payoffs.

In this context, when Monty plays $\QQ$, Conies' Bayesian strategy $\PP'=\PP'(\QQ)$ is called a {\it best response}
to $\QQ$. Reciprocally, 
when  Conie plays a particular $\PP$ Monte's best response strategy $\QQ'=\QQ'(\PP)$
is the one for which Conie's winning probability is minimal,
$$\PP\CC\QQ'^T= \min_{\PP}  \PP\CC\QQ^T.$$

A profile of actors' mixed strategies $(\PP^*,\QQ^*)$ is  said to be a {\it minimax solution}  (or saddle point)
if the  strategies are the best responses to each other,
$$ \PP^*\CC\QQ^{*T}= \max_{\PP}  \PP\CC\QQ^{*T} = \min_{\QQ} \PP^*\CC\QQ^T.$$
Such a solution exists for arbitrary matrix game with finite sets of strategies, by the minimax theorem due to von Neumann.
The number 
$V:=\PP^*\CC\QQ^{*T}$ is called the {\it  value of the game}, and it does not depend on particular minimax profile due to the fundamental
relation
$$V=\max_{\PP}\min_{\QQ}  \PP\CC\QQ^{T} = \min_{\QQ}\max_{\PP} \PP\CC\QQ^T$$
involved in the minimax theorem.

Strategy $\PP^*$ of Conie  is minimax if it guarantees the winning probability at least $V$ no matter what Monte does.
Strategy $\QQ^*$ of Monte is minimax if it incures the winning probability at most $V$ no matter what Conie does.

Recalling our discussion around the dominance it is really easy to see that $V=2/3$.
If Monte chooses $\theta$ uniformly at random, Conie cannot winning probability higher than $2/3$. On the other hand,
if $x$ is chosen uniformly at random and $x\ts\ts$ is played, Conie guarantees $2/3$ no matter what Monte does.
So a solution is found.

We wish to approach the solution more formally, manipulaing the payoff matrix.
The principle of the elimination of dominated strategies says:

\begin{itemize}
\item[] the value of the game is not altered if the game matrix is (repeatedly) reduced by eliminating a (weakly) dominated row or column.
\end{itemize}

Explicitly, strategy $1\swsw$~ dominates 2\masw~ and 2\mama~  
\begin{center}
\begin{tabular}{c|cccccc}
1\swsw &0  & 0 &  1  &  1 & 1  &1\\
2\masw &0  & 0 &  1 &  0 & 1  &1\\
2\mama &0  & 0 &  1  &  1 & 0  &0\\
\end{tabular}
\end{center}
and 3\swsw~ dominates 2\swma~
\begin{center}
\begin{tabular}{c|cccccc}
3\swsw &1  & 1 &  1  &  1 & 0  &0\\
2\swma &1  & 1 &  0  &  1 & 0  &0\\
\end{tabular}
\end{center}
Continuing the elimination we reduce to
\begin{center}
\begin{tabular}{c|cccccc}
1\swsw &0  & 0 &  1  &  1 & 1  &1\\
2\swsw &1  & 1 &  0  &  0 & 1  &1\\
3\swsw &1  & 1 &  1  &  1 & 0  &0\\
\end{tabular}
\end{center}
Finally, discarding repeated columns
the matrix is reduced to the square matrix $\Cc$:
\vskip0.3cm
\begin{center}

\begin{tabular}{c|ccc}
$\theta,y=$     & 12  & 21 & 32 \\
\hline
1\swsw &0  &   1   & 1  \\
2\swsw &1  &   0   & 1  \\
3\swsw &1  &   1   & 0  \\
\end{tabular}
\end{center}
\vskip0.3cm
The reduced game $\Cc$ has a clear interpretation.
Monte and Conie choose door numbers $\theta$ and $x$, respectively, from $\{1,2,3\}$. 
If the choices {\it mismatch} ($\theta\neq x$) Conie wins, otherwise there is no payoff.

Let us find an {\it equalizing} strategy $\Pp^*$ for which ${\Pp}^*\Cc{\Qq}^T$ is the same no matter 
which  counter-strategy $\Qq$ Monte plays. Taking for $\Qq$ pure strategies we  arrive at the system of equations
$$p_2+p_3=p_1+p_2=p_1+p_3,$$
which taken together with $p_1+p_2+p_3=1$ is complete. 
Solving the system we see that with $\Pp^*=(1/3,1/3,1/3)$ Conie wins with probability 2/3
no matter what Monte does. Similarly, when Monte plays $\Qq^*=(1/3,1/3,1/3)$ the winning probability is always 
2/3 no matter what Conie does. Thus $(\Pp^*,\Qq^*)$ is a solution of the game and the value $V=2/3$ is confirmed.

Yet another way to arrive at $(\Pp^*,\Qq^*)$ is to use symmetry of matrix $\Cc$ induced by permutations of the door numbers
(see \cite{Ferguson}, Theorem 3.4).

A related game with diagonal matrix 
$$\left(
\begin{array}{ccc}
-1  &   0   & 0  \\
 0  &   -1   & 0  \\
 0  &   0   & -1  \\
\end{array}\right)
$$
is obtained by subtracting a constant matrix from $\Cc$. In this diagonal game {\it Monte} wins if the choices match, otherwise there is no payoff.

Going back to the original matrix $\CC$, we conclude 
that the profile 
\begin{eqnarray*}
\PP^*=\left({1\over 3},0,0,0,{1\over 3},0,0,0,{1\over 3},0,0,0\right),~~
 \QQ^*_{1,1,1}=\left({1\over 3},0,{1\over 3},0,{1\over 3},0\right)
\end{eqnarray*}
is a solution to the game. Strategy $\PP^*$ is equalizing, as $\PP^*\CC \QQ$ =2/3 for all $\QQ$.
According to the solution $(\PP^*, \QQ^*_{1,1,1})$ Monte plays the `crawl' strategy: he
hides the prize uniformly at random,
and he always reveals the door with higher number,  when there is a freedom for the second action. 
Conie selects door $x$ uniformly at random and always plays $\ts$.



The generic Monte's strategy assigning  probability $1/3$ to every value of $\theta$ is of the form 
$$\QQ^*_{{\lambda_1},{\lambda_2},{\lambda_3}}=\left({\lambda_1\over 3}\,,\,{1-\lambda_1\over 3}\,,\,
{\lambda_2\over 3}\,,\,{1-\lambda_2\over 3}\,,\,{\lambda_3\over 3}\,,\,{1-\lambda_3\over 3}\right)$$
where $0\leq\lambda_\theta\leq 1$, $\theta\in \{1,2,3\}$.
Parameter $\lambda_\theta$ is the conditional 
probability that Monte will offer switching to door with smaller number in the case of match. 
If Conie plays best response to $\QQ^*_{{\lambda_1},{\lambda_2},{\lambda_3}}$ she wins with probability $2/3$, as 
we have seen when considering the Bayesian setting. Therefore  each strategy $\QQ^*_{{\lambda_1},{\lambda_2},{\lambda_3}}$ is minimax.
On the other hand, if the values of $\theta$ have probabilities $\pi_\theta$ a best response of Conie yields winning probability
$1-\min(\pi_1,\pi_2,\pi_3)$, which is
minimized for $\pi_1=\pi_2=\pi_3=1/3$. We see that
\begin{itemize}
\item[] The strategies $\QQ^*_{{\lambda_1},{\lambda_2},{\lambda_3}}$ and only they are the minimax strategies of Monte.
\end{itemize}

Furthermore, if $\QQ^*_{{\lambda_1},{\lambda_2},{\lambda_3}}$ is fully supported, which is the case when
$\lambda_{\theta}\in (0,1)$ for $\theta\in \{1,2,3\}$,
then the unique best response is $\PP^*$, thus
\begin{itemize}
\item[] Conie's strategy $\PP^*$ of choosing a door uniformly at random, then always switching is the unique minimax strategy.
\end{itemize}

The subclass of Monte's strategies with the second action independent of the first given match consists of strategies with equal probabilities
$\lambda_1=\lambda_2=\lambda_3$; these were discussed e.g. in  \cite{Rosenhouse} (Version Five of the MHP).

We note that the general theory does not preclude 
weakly dominated strategies from being minimax
(see \cite{Ferguson}, Section 2.6, Exercise 9). This does not occur in the MHP because there exist fully supported minimax strategies of Monte.

\section{The general-sum games}

The logical way to  go beyond the zero-sum game is a non-zero sum game.
In such model there is a $12\times 6$ payoff matrix of Monte $\HH$ which need not be the negative of $\CC$.
A best response to Monte's strategy $\QQ$ is defined as before, but best response to Conie's strategy $\PP$ is now
a strategy which maximizes the expected payoff $\PP\HH\QQ^T$.

A  central solution concept for the general-sum game game 
is a {\it Nash equilibrium}, defined as 
a profile of mixed strategies $(\PP',\QQ')$ which are best responses to each other,
$$\PP'\CC\QQ'^{T}=\max_{\PP} \PP\CC\QQ'^{T} ~~~{\rm  and~~~} \PP'\HH\QQ'^{T}=\max_{\QQ} \PP'\HH\QQ^{T}.$$
That is to say, Nash equilibrium is a profile 
 $(\PP',\QQ')$ such that a unilateral change of the strategy by one of the actors cannot improve the private 
payoff of the actor.
A general theorem due to John Nash ensures that at least one such Nash equilibrium exists.

Nash equilibrium is a concept of {\it noncooperative} game theory. The players cannot make binding agreements on a joint choice of strategy unless 
the agreements are self-enforced. This is formalized by the independence presumed in the product formulas for computing the payoffs.

In every Nash equilibrium Conie will have  winning probability not less than her minimax value $V=2/3$, which is 
her {\it safety level}. Higher probability in some Nash equilibria might be possible, since the game is no longer  antagonistic.

Our  analysis of the Bayesian strategies can be applied to the general-sum games as well. Suppose in a Nash equilibrium $(\PP',\QQ')$ 
strategy $\QQ'$ is fully supported. Then $\PP'$, being a best response to $\QQ'$, is a mixture of always-switching strategies.
Let $\pi_1\geq \pi_2\geq \pi_3>0$ be probabilities of the values $\theta=1,2,3$ under $\QQ'$.
Then we have

\begin{itemize}
\item[]
{\rm A profile $(\PP',\QQ')$ is a Nash equilibrium with fully supported  
$\QQ'$ if and only if there exists a probability vector $(p_1,p_2,p_3)$
 such that  the mixture  with weights $p_1, p_2, p_3$  
of  rows $1\ts\ts, 2\ts\ts, 3\ts\ts$  of matrix $\HH$ is a row vector with equal entries.  There are three 
possibilities
\vskip0.1cm
\noindent
\begin{itemize}
\item[\rm(1)]  $p_3=1$, $\pi_2>\pi_3$ and the row $3\ts\ts$ of $\HH$ is a constant row,  
\item[\rm(2)]  $p_1=0$, $\pi_1>  \pi_2=\pi_3$ and a mixture of the rows $2\ts\ts$ and $3\ts\ts$ of $\HH$ is a constant row,
\item[\rm(3)]  $p_1p_2p_3>0,~ \pi_1=\pi_2=\pi_3=1/3$, and the arithmetical average of rows 
$1\ts\ts$, $2\ts\ts$ and $3\ts\ts$ of $\HH$ is a constant row
\end{itemize}
}
\end{itemize}
In the case (3)
Conie's winning  probability is her safety level 2/3 but, unlike the zero-sum game,  
the equilibrium strategy $\PP'$ need not give the same  probability to every always-switching strategy.
Thus if the equilibrium has a property of nondegeneracy, when every pure strategy of Monte and every always-switching strategy of Conie
is played with positive probability, then for Conie the game brings the same as the game against antagonistic Monte.

If none of the mixtures of the $x\ts\ts$-rows is constant, there is no fully supported Nash equilibrium.

What could be plausible assumptions on Monte's payoff $\HH$? If Monte is only concerned about the fate of the prize,
and not where and how the prize is won,
 his essentially distinct payoff structures are
$\HH=-\CC$ (antagonistic Monte) $\HH=\CC$ (sympathetic Monte) and $\HH={\bf 0}$ (indifferent Monte). The first case is zero-sum,
in the second case every entry 1 of $\CC$ corresponds to a pure Nash equilibrium, and in the third case every pair 
$(\QQ, \PP')$ with best-response $\PP'=\PP'(\QQ)$ is a Nash equilibrium.

It is not hard to design further exotic examples of payoffs $\HH$ for which context-dependent strategies  enter some Nash equilibrium profile.

If there is a moral of that it is perhaps this: context-dependent strategies 
{\it may} be a rational kind of behavior under certain intensions of Monte.
 The
source for this phenomenon is twofold. Firstly, 
this is the very idea of equilibrium:
if Conie steps away from the context-dependent equilibrium strategy, 
Monty is not forced to play old strategy and may change biases in a way unfavorable for Conie, typically pushing her down to the safety level $2/3$.
Secondly, the domination is only {\it weak}, thus can have no effect if support is not full.
However, such context-dependent equilibria are highly unstable, and minor perturbations of $\HH$ will destroy them.
Practically speaking, if Conie has any doubts about Monte's intensions 
it is safe to stay with always-switching strategies.

As for the famous question, the noncooperative game theory gives more weight to vos Savant's solution by adding

\begin{itemize}
\item[]
{\it Yes, you should switch. You knew the rules of the game.
If your decision were to pick door 1 and hold when a switch to door 2 offered, 
then I could beat your strategy by picking door 2 and switching whichever happens.
My strategy will be even strictly better than yours if  the prize can ever be hidden behind door 2.}
\end{itemize}

If the game is {\it cooperative}, for instance if Monte and Conie want happily drive their new cadillac to Nice,
they could just favor door 1. But this is a completely different story.


\begin{thebibliography}{99}

\bibitem{Dekking} Dekking, F.M. Kraaikamp, C.,  Meester, L.E. and Lopuha{\"a}, H.P.
{\it A modern introduction in probability and statistics: understanding why and how}, Springer, 2005. 


\bibitem{Ferguson}Ferguson T. (2000) {\it A course on game theory}.\\ {\tt  http://www.math.ucla.edu}



\bibitem{Doors} Gnedin, A. (2011) The doors, {\tt arxiv.org}

\bibitem{Dominance} Gnedin, A. (2011) Dominance in the Monty Hall problem. {\tt arxiv.org}


\bibitem{Ethier} Ethier, S.N. {\it The doctrine of chances: probabilistic aspects of gambling}, Springer, 2010.


\bibitem{Gill} Gill, R. (2011) The Monty Hall problem is not
a probability puzzle (it's a challenge in mathematical modelling), 
{\it Statistica Neerlandica} {\bf 65} 58--71.

\bibitem{Grinstead} Grinstead, C.M. and Snell, J. L. 
{\it Grinstead and Snell's introduction to probability}. \\
{\tt http://www.math.dartmouth.edu/~prob/prob/prob.pdf}  2006
(online version of Introduction to Probability, 
2nd edition, published by the Amer. Math. Soc.).


\bibitem{Haggstrom} H\"aggstr\"om O. 
{\it Streifz\"uge durch die Wahrscheinlichkeitstheorie}, Springer 2005. 


\bibitem{Rosenhouse}Rosenhouse, J. 
{\it The  Monty Hall problem}, Oxford University Press, 2009.
\bibitem{Rosenthal} 
Rosenthal, J.S. (2008) 
Monty Hall, Monty Fall, Monty Crawl. {\it Math Horizons} 5-7, September issue.
\bibitem{Tijms} Tijms, H. {\it Understanding probability}, CUP, 2007.
\end{thebibliography}
\end{document}